# Multi-Criteria Shape Optimization of Flow Fields for Electrochemical Cells

Dr. Sebastian Blauth[1,*], Dr. Marco Baldan[1], Dr. Sebastian Osterroth[1], Dr. Christian Leithäuser[1], Prof. Dr. Ulf-Peter Apfel[2,3], Julian Kleinhaus[3], Kevinjeorjios Pellumbi[2,3], Dr. Daniel Siegmund[2,3], Dr. Konrad Steiner[1], Prof. Dr. Michael Bortz[1]

We consider the shape optimization of flow fields for electrochemical cells. Our goal is to improve the cell by modifying the shape of its flow field. To do so, we introduce simulation models of the flow field with and without the porous transport layer. The latter is less detailed and used for shape optimization, whereas the former is used to validate our obtained results. We propose three objective functions based on the uniformity of the flow and residence time as well as the wall shear stress. After considering the respective optimization problems separately, we use techniques from multi-criteria optimization to treat the conflicting objective functions systematically. Our results highlight the potential of our approach for generating novel flow field designs for electrochemical cells.

**Keywords:**

CFD; Electrochemical Cell; Multi-Criteria Optimization; Numerical Optimization; Shape Optimization

**Author affiliations**

[1]Fraunhofer ITWM, Fraunhofer-Platz 1, 67659 Kaiserslautern, Germany.

[2]Fraunhofer UMSICHT, Osterfelder Str. 3, 46047 Oberhausen, Germany.

[3]Ruhr-Universität Bochum, Universitätsstr. 150, 44801 Bochum, Germany.

Email corresponding author: sebastian.blauth@itwm.fraunhofer.de

ORCID iDs of the authors

0000-0002-1577-2420 (Ulf-Peter Apfel)

0000-0002-5803-3150 (Marco Baldan)

0000-0001-9173-0866 (Sebastian Blauth)

0000-0001-8169-2907 (Michael Bortz)

0000-0002-0901-4789 (Julian Kleinhaus)

0000-0001-8936-9805 (Christian Leithäuser)

0000-0001-9702-521X (Sebastian Osterroth)

0000-0002-5313-5799 (Kevinjeorjios Pellumbi)

0000-0003-2476-8965 (Daniel Siegmund)

# 1  Introduction

To achieve the large-scale establishment of electrolysis technologies, the optimization of electrolyzers and their components is essential [1]. Within these technologies, hydrogenation reactions play a central role. Electrochemical hydrogenation reaction (EChH) emerged in recent years as a more sustainable alternative to the current thermocatalytic state of the art [2–9]. Instead of operating at elevated temperatures and pressures, EChH electrolyzers operate at ambient conditions using water and, ideally, renewable electricity to selectively hydrogenate unsaturated bonds. We have recently demonstrated the efficient production of the vitamin synthon 2-methyl-3-buten-2-ol in scalable zero-gap electrolyzers (ZGEs), showing how material selection plays a critical role in the observed activity [10, 11]. Compared to other electrolytic cell types, ZGEs use a solid polymer electrolyte to conduct ions between the two half-cell compartments and minimize ohmic resistance by mechanically compressing the two electrodes on the polymer electrolyte [12, 13]. Despite being one of the most promising cell architectures for EChH, there is currently little information on how to improve the efficiency and design of components of ZGEs for electrosynthetic applications [10, 11, 14, 15].

In this work, we investigate how methods from mathematical modeling, simulation, and optimization can be applied to improve the performance of EChH reactors. We focus on the task of optimizing flow field geometries for the cathode side of the EChH reactor. The methods used to achieve this task are computational fluid dynamics (CFD), shape optimization with partial differential equations (PDEs) and multi-criteria optimization (MCO). The goal of shape optimization is to improve the quality of some system by altering its shape, in this case, the shape of the flow field. We employ techniques from shape calculus to perform shape optimization efficiently without requiring parametrizations. Such approaches have been applied in the literature to optimize aircraft [16], electromagnetic devices [17, 18], polymer spin packs [19, 20], and microchannel systems [21, 22]. Additionally, there are many recent developments in the field of shape optimization which enable the optimization of complex industrial applications. This can be seen, e.g., in [23], where space mapping methods for shape optimization are introduced, in [24], where the shape optimization software package cashocs is presented, in [25, 26], where efficient solution algorithms for shape optimization are investigated, or in [27], where sophisticated mesh deformations for shape optimization are proposed.

Multi-criteria optimization in the context of shape optimization has already received some attention in the literature, e.g., in [28–31]. In these publications, the shapes of the boundaries are parametrized using, e.g., B-splines, which limits the set of reachable designs. However, in this work, we consider free-form shape optimization without parametrizations of the geometry.

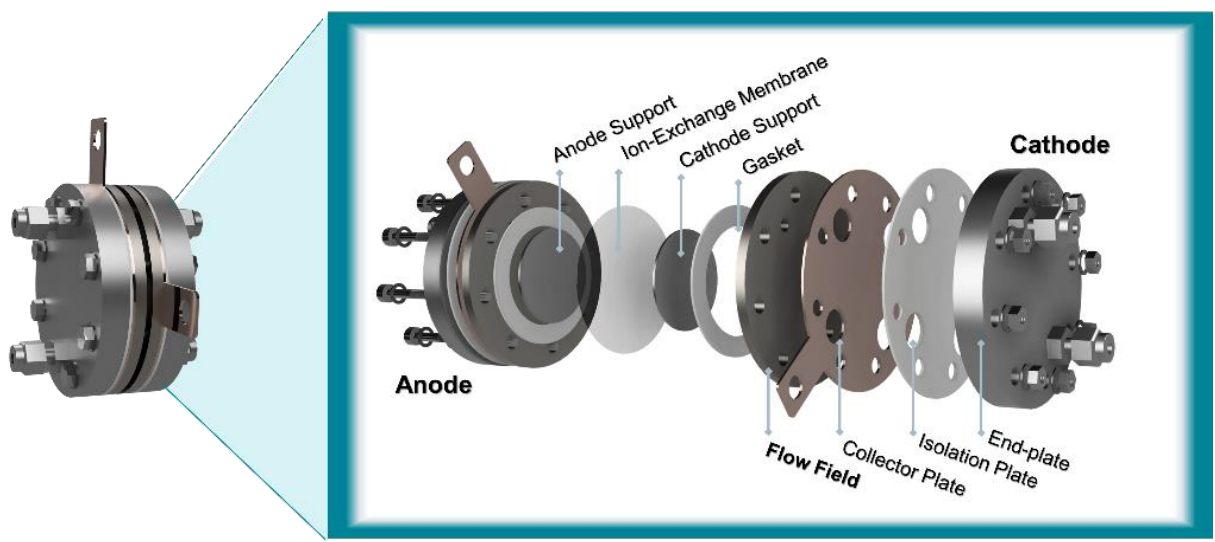

**Figure 1.** Assembly of an EChH reactor.

## 2 Methods and Models

The assembly of an EChH reactor is visualized in Figure 1. We focus our study on the cathode side of the reactor. Water is fed into the flow field geometry and some of it enters the porous transport layer (PTL, marked as cathode support in Figure 1). Then the electrochemical reaction takes place and ions are exchanged through the ion-exchange membrane. Flow field and PTL are modeled in approximation of our recently reported EChH setup with a graphite felt PTL (Sigracell GFD 2.5 EA, SGL Carbon) [10].

We first introduce a model of the flow through the combination of flow field and PTL based on CFD methods. To optimize the shapes of flow fields, we use a simplified model which only considers the flow field without the PTL. Based on this model, we perform a shape optimization with different objectives to derive new flow field shapes. Methods from multi-criteria optimization are used to investigate the combination of multiple conflicting objectives. Finally, the flow field shapes are evaluated with the detailed simulation model.

### 2.1 Flow Field and PTL Simulation Model

The flow through a porous medium can, in general, be modeled on different scales. When the geometric structure of the porous medium is resolved, this is called a detailed or microscale model. However, the use of such a model is expensive. Therefore, the porous medium can be described as effective medium by using a permeability tensor, resulting in a so-called macroscopic model. As we consider flow in both porous and pure fluid domains simultaneously, we employ the commonly used Brinkman model [32].

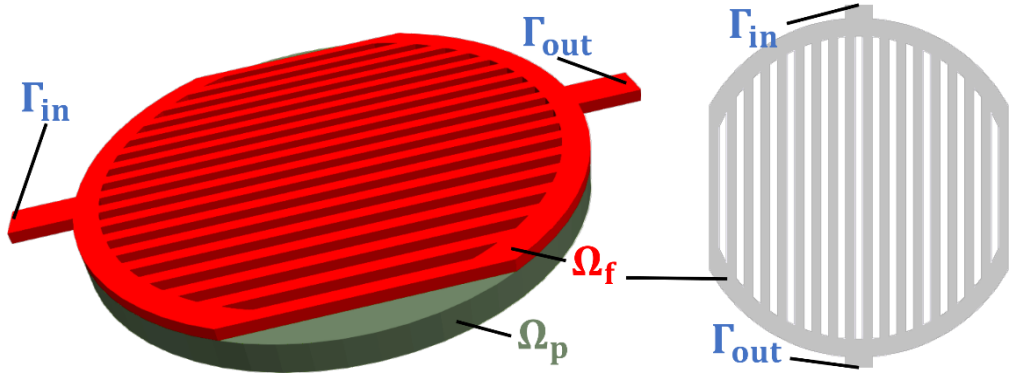

**Figure 2.** Geometry for the detailed (left) and simplified (right) flow field simulation model. This is the original flow field shape with parallel channels which is denoted by PAR throughout this manuscript.

The geometric setting for the detailed simulation model is shown on the left in Figure 2. The computational domain $\Omega$ is split into a pure fluid part $\Omega_\text{f}$ and a porous part $\Omega_\text{p}$. In the following, we restrict our considerations to steady-state, laminar, isothermal, and incompressible flow.

The flow in $\Omega_f$ is modeled using the Navier-Stokes equations. For $\Omega_p$, a viscous resistance term is added, which results from Darcy's law. The equations read

$$\begin{aligned}
\nabla \cdot u &= 0 \quad \text{in } \Omega, \\
\rho(u \cdot \nabla)u - \nabla \cdot (\mu \nabla u) + \nabla p &= 0 \quad \text{in } \Omega_\text{f}, \\
\rho(u \cdot \nabla)u - \nabla \cdot (\mu_\text{eff} \nabla u) + \mu K^{-1} u + \nabla p &= 0 \quad \text{in } \Omega_\text{p}.
\end{aligned} \quad (1)$$

Here $u$ denotes the velocity and $p$ the pressure. The properties $\rho$ and $\mu$ describe the density and dynamic viscosity of the fluid, respectively. In the porous part, the viscosity might be different and, therefore, the effective viscosity $\mu_\text{eff}$ is used. However, there is no consensus about the definition of the effective viscosity so that mostly the value of $\mu$ is used [33]. The porous medium is characterized by its permeability $K$.

The boundary of $\Omega$ is divided into the inlet $\Gamma_{\text{in}}$, the wall boundary $\Gamma_{\text{wall}}$, and the outlet $\Gamma_{\text{out}}$. Due to the symmetry of the geometry, we only use half of it and denote the corresponding part of the boundary as $\Gamma_{\text{sym}}$. Moreover, we denote by $\Gamma_{\text{int}}$ the interface between the fluid and porous part.

The following boundary conditions are specified:

$$\begin{aligned} u &= u_{\text{in}} \quad \text{on } \Gamma_{\text{in}}, \\ u &= 0 \quad \text{on } \Gamma_{\text{wall}}, \\ \mu\, \partial_n u - p n &= 0 \quad \text{on } \Gamma_{\text{out}}, \\ u \cdot n &= 0 \quad \text{on } \Gamma_{\text{sym}}, \\ \mu\, \partial_n u \times n &= 0 \quad \text{on } \Gamma_{\text{sym}}. \end{aligned} \quad (2)$$

Here, $\partial_n u$ denotes the normal derivative of a vector field $u$, which is defined as $\partial_n u = Du\, n$, i.e., the application of the Jacobian $Du$ to $n$, where the latter is the unit outer normal vector on $\Gamma$.

For the computation of velocity and the pressure, the software FiltEST [34] is used.

## 2.2 Flow Field Simulation Model for Shape Optimization

For the shape optimization of the flow field, we restrict our attention to the fluid domain $\Omega_f$ only and, for the moment, neglect the porous part $\Omega_p$. Hence, we employ a no-slip boundary condition on the interface $\Gamma_{\text{int}}$ between fluid and porous part. Our simplified model for the flow in the flow field is given by

$$\begin{aligned} \rho(u \cdot \nabla)u - \mu \Delta u + \nabla p &= 0 \quad \text{in } \Omega_f, \\ \nabla \cdot u &= 0 \quad \text{in } \Omega_f, \\ u &= u_{\text{in}} \quad \text{on } \Gamma_{\text{in}}, \\ u &= 0 \quad \text{on } \Gamma_{\text{wall}} \cup \Gamma_{\text{int}}, \\ u \cdot n &= 0 \quad \text{on } \Gamma_{\text{sym}}, \\ \mu\, \partial_n u \times n &= 0 \quad \text{on } \Gamma_{\text{sym}}, \\ \mu\, \partial_n u - p n &= 0 \quad \text{on } \Gamma_{\text{out}}. \end{aligned} \quad (3)$$

For the sake of faster simulation times, we solve the above model with the help of a dimension reduction technique presented in [21], which reduces the above three-dimensional problem to a two-dimensional one, significantly reducing the cost of its numerical solution while sacrificing only little accuracy. We refer the reader to [21] for more details and a thorough investigation of this dimension reduction technique.

## 2.3 Cost Functions for the Shape Optimization

Based on the simplified flow model (3) we optimize the shape of the flow field. To do so, we introduce three cost functions $J_1, J_2$, and $J_3$ with the following goals:

- $J_1$ aims at achieving a uniform flow distribution among the channels of the flow field,
- $J_2$ aims at achieving a uniform residence time of the flow in each channel, and
- $J_3$ aims to control the shape of the in- and outflow regions via a criterion based on the wall shear stress.

### 2.3.1 Uniform Flow Cost Function $J_1$

Motivated by the observation that the volume flow rate varies significantly between the channels of the original domain PAR (cf. Section 3.2), we define the cost function

$$J_1(\Omega_f, u) = \tfrac{1}{2} \int_{\Omega_{\text{channels}}} (u - u_{\text{des}})^2 \mathrm{d}x. \quad (4)$$

Here, $\Omega_{\text{channels}}$ is the part of the domain $\Omega_f$ corresponding to the channels and $u_{\text{des}}$ is the desired velocity, which we want to achieve within the channels. The latter can be computed a-priori as Poiseuille flow profile for a single channel with the desired average volume flow rate $\dot{V}_{des}$ (cf. [21]). The cost functional aims to achieve a uniform flow distribution among the channels.

### 2.3.2   Uniform Residence Time Cost Function $J_2$

The second cost function $J_2$ is given by

$$J_2(\Omega_f, u) = \frac{1}{2}\sum_i (\tau_i - \tau_{des})^2, \tag{5}$$

where $\tau_i$ is the residence time in channel $i$ and $\tau_{des}$ is a desired residence time. Hence, the cost functional aims at achieving a uniform residence time, $\tau_{des}$, in each channel. In this paper we restrict ourselves to the case of a uniform residence time for all channels, but using different desired residence times for different channels would be possible. The residence time for channel $i$ is computed as

$$\tau_i = \frac{L_i}{\bar{v}_i}, \tag{6}$$

with $L_i$ and $\bar{v}_i$ being the length and average velocity of channel $i$, respectively, with

$$\bar{v}_i = \frac{\left|\int_{\Omega_i} u_y dx\right|}{\int_{\Omega_i} 1 dx}. \tag{7}$$

Here, $u_y$ is the velocity component in channel direction and $\Omega_i$ is the subdomain corresponding to channel $i$.

### 2.3.3   Wall Shear Stress Cost Function $J_3$

As third objective we consider the wall shear stress on the outer boundaries of the flow field. The goal of this criterion is to prevent the creation of unnecessarily large inflow and outflow volumes during the optimization. This is achieved by controlling the wall shear stress on the boundary part $\Gamma_{wss} = \Gamma_{wss}^{in} \cup \Gamma_{wss}^{out}$ (cf. Figure 5). The wall shear stress can be seen as a measure of how fast the flow is near a boundary, hence, a low wall shear stress can indicate a stagnation zone. So, a criterion that ensures a sufficiently high wall shear stress can prevent stagnation zones and, thus, unnecessarily large volumes. For this reason, we introduce the constraint

$$\sigma \geq \sigma_{thr} \quad \text{on } \Gamma_{wss}, \tag{8}$$

where $\sigma$ is the wall shear stress, which is defined as

$$\sigma = \mu \|\partial_n u\|. \tag{9}$$

Moreover, $\sigma_{thr}$ is a threshold value that should be reached on $\Gamma_{wss}$. To treat this constraint numerically, we employ a Moreau-Yosida regularization (cf. [35]) and end up with the cost functional

$$J_3(\Omega_f, u) = \int_{\Gamma_{wss}} \min(0, \sigma - \sigma_{thr})^2 \, ds. \tag{10}$$

For more details, we refer the reader, e.g., to [20], where a similar criterion is used in the context of designing polymer spin packs without stagnation zones.

### 2.4   Numerical Solution of the Shape Optimization Problems

First, we consider the single criteria shape optimization problem

$$\min_{\Omega_f, u} J_i(\Omega_f, u) \tag{11}$$
$$\text{subject to (3),}$$

for $J_i$ being any of the previously defined cost functions $J_1$, $J_2$ or $J_3$.

To solve the optimization problem (11) numerically, we use our shape optimization software cashocs [24], which is a software package for the automated solution of PDE-constrained shape optimization problems. We impose the following geometrical constraints on the problem: The outer boundary of the geometry, including inlet, outlet, and symmetry boundary, is not allowed to change. Additionally, the channels of the geometry are only allowed to change in length. As cashocs is based on the finite element software FEniCS [36], we use the finite element method to discretize the state and adjoint systems with LBB-stable Taylor-Hood elements, i.e., piecewise quadratic and linear Lagrange elements for velocity and pressure, respectively. For the solution of the optimization problem, we use a limited memory BFGS method from [25] which is implemented in cashocs.

## 2.5 Multi-Criteria Shape Optimization

Finally, we can formulate the PDE-constrained multi-criteria shape optimization problem as follows:

$$\min_{\Omega_f, u} J(\Omega_f, u) = [J_1(\Omega_f, u), J_2(\Omega_f, u), J_3(\Omega_f, u)] \tag{12}$$

subject to (3),

where we consider the above multi-criteria optimization problem (12) in the sense of Pareto-efficiency (cf. [37]).

To approximate the Pareto front, the calculation of single Pareto points is performed by scalarizing the cost functions to arrive at a single objective function to be minimized. To do so, the weighted-sum scalarization is used [37]:

$$\min_{\Omega_f, u} J_\lambda(\Omega_f, u) = \frac{\lambda_1}{J_1(\Omega_{f,0}, u_0)} J_1(\Omega_f, u) + \frac{\lambda_2}{J_2(\Omega_{f,0}, u_0)} J_2(\Omega_f, u) + \frac{\lambda_3}{J_3(\Omega_{f,0}, u_0)} J_3(\Omega_f, u) \tag{13}$$

subject to (3),

where $\lambda_i$ are non-negative weights. To ensure that the problem is well-scaled, each cost functional is divided by its respective value for the initial domain $\Omega_{f,0}$ with corresponding velocity $u_0$, where the original flow field shape PAR is used as $\Omega_{f,0}$. An effective approximation of the Pareto front is obtained by applying the so-called sandwiching algorithm [38, 39]. Sandwiching relies on an iterative approach that provides, at each step, a weight vector $\lambda = [\lambda_1, \lambda_2, \lambda_3]$ to be used in (13). New Pareto points are added as long as the desired approximation quality [39] is not achieved. For the sake of brevity, we do not give the details of the sandwiching technique here but refer the reader to [38, 39].

The decision maker can then interactively explore the Pareto set by navigating with a graphical user interface [38]. Graphical sliders are used for this purpose, each corresponding to one cost function. By moving one slider, the other sliders are updated in real time, i.e., information on the trade-offs between the best compromises is directly visualized. We refer the reader to, e.g., [38–40] for more details.

## 3 Results and Discussion

We begin our study by investigating the original flow field shape PAR using our detailed simulation model. We then use shape optimization and multi-criteria shape optimization to derive optimized flow field shapes which are compared based on criteria derived from the detailed simulation model.

### 3.1 Simulation Setup

For now, we only consider the flow of water through the flow field and PTL and neglect the electrochemical reaction. Note that, in reality, the electrochemical reaction produces gas as a side product, which mixes with the water and thus changes the flow properties when advancing through the cell. It is up for future work to potentially include this effect into the simulation model.

The fluid is modeled using a constant viscosity $\mu = 3.547 \cdot 10^{-4}$ Pa·s and a constant density $\rho = 971.79 \frac{kg}{m^3}$, which corresponds to water at 80°C. The inlet velocity is chosen so that we have an inlet flow rate of $\dot{V}_{in} = 7.5 \frac{mL}{min}$. The permeability of the PTL is estimated as $K = 10^{-11} \ m^2$ from [41].

## 3.2 Simulation Results for Original Flow Field PAR

Simulation results obtained from the detailed simulation model (1) are shown in Figure 3. The left plot shows streamlines in the flow field and the right plot shows streamlines in the PTL. Both plots are colored by the velocity magnitude. Note that in the plot on the right, there is an area with few streamlines. This area is still covered by flow, but only reached by a few of the streamlines used for the visualization.

These simulations show an uneven flow distribution: while the flow rate in the middle two channels is high, the remaining channels receive significantly less flow. This also translates to the flow in the PTL, which is also much higher in the middle. Based on this observation, we focus on the improvement of flow uniformity in the following. In addition, we also investigate the influence of other criteria, such as residence time and wall shear stress.

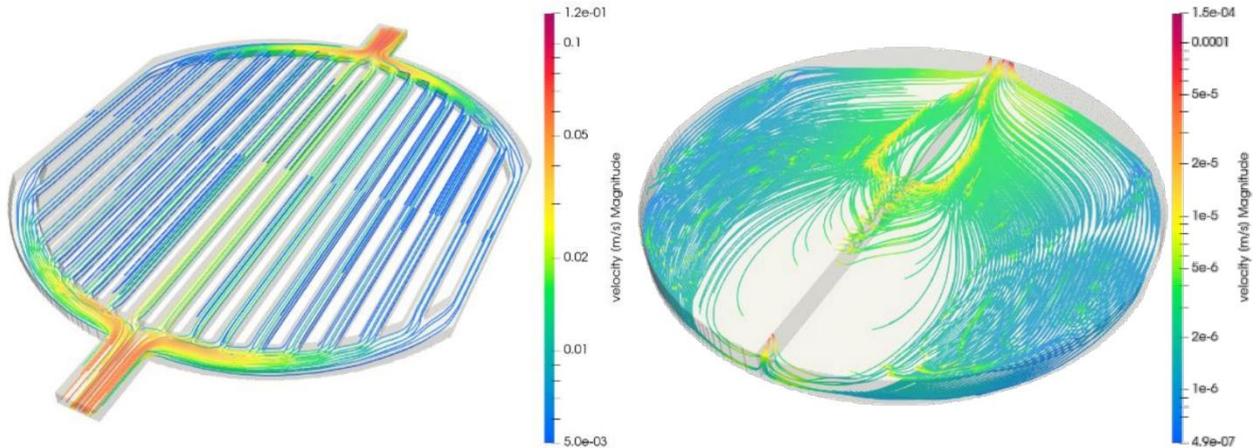

**Figure 3.** Simulation of the original flow field PAR using the detailed simulation model. Streamline visualization of velocity magnitude for flow field (left) and PTL (right).

## 3.3 Shape Optimization of Flow Field

To improve the flow properties of the flow field, we solve the single criteria shape optimization problem (11) for each of the previously defined cost functions: we use the uniform flow cost function $J_1$ but also test the influence of the two other cost functions $J_2$ and $J_3$. Solving the shape optimization problem gives rise to the following shapes:

- **PAR J1**: Uniform flow per channel, minimizing the cost function $J_1$,
- **PAR J2**: Uniform residence time per channel, minimizing the cost function $J_2$, and
- **PAR J3**: Sufficiently high wall shear stress on $\Gamma_{wss}$, minimizing the cost function $J_3$.

The resulting flow field shapes are shown in Figure 4.

Figure 5 compares the individual quantities which give rise to the cost functions:

1) The first column in Figure 5 shows the volume flow rate through each of the 18 channels of the flow field, corresponding to cost function $J_1$. The volume flow rate is obtained from the velocity integrated over the cross section of each channel.
2) The second column in Figure 5 shows the residence time within each of the 18 channels of the flow field, corresponding to cost function $J_2$.

3) The third column in Figure 5 shows the wall shear stress on the upper wall $\Gamma_{wss}^{in}$ of the flow field, which contributes to cost function $J_3$.
4) The fourth column in Figure 5 shows the wall shear stress on the lower wall $\Gamma_{wss}^{out}$ of the flow field, which contributes to cost function $J_3$.

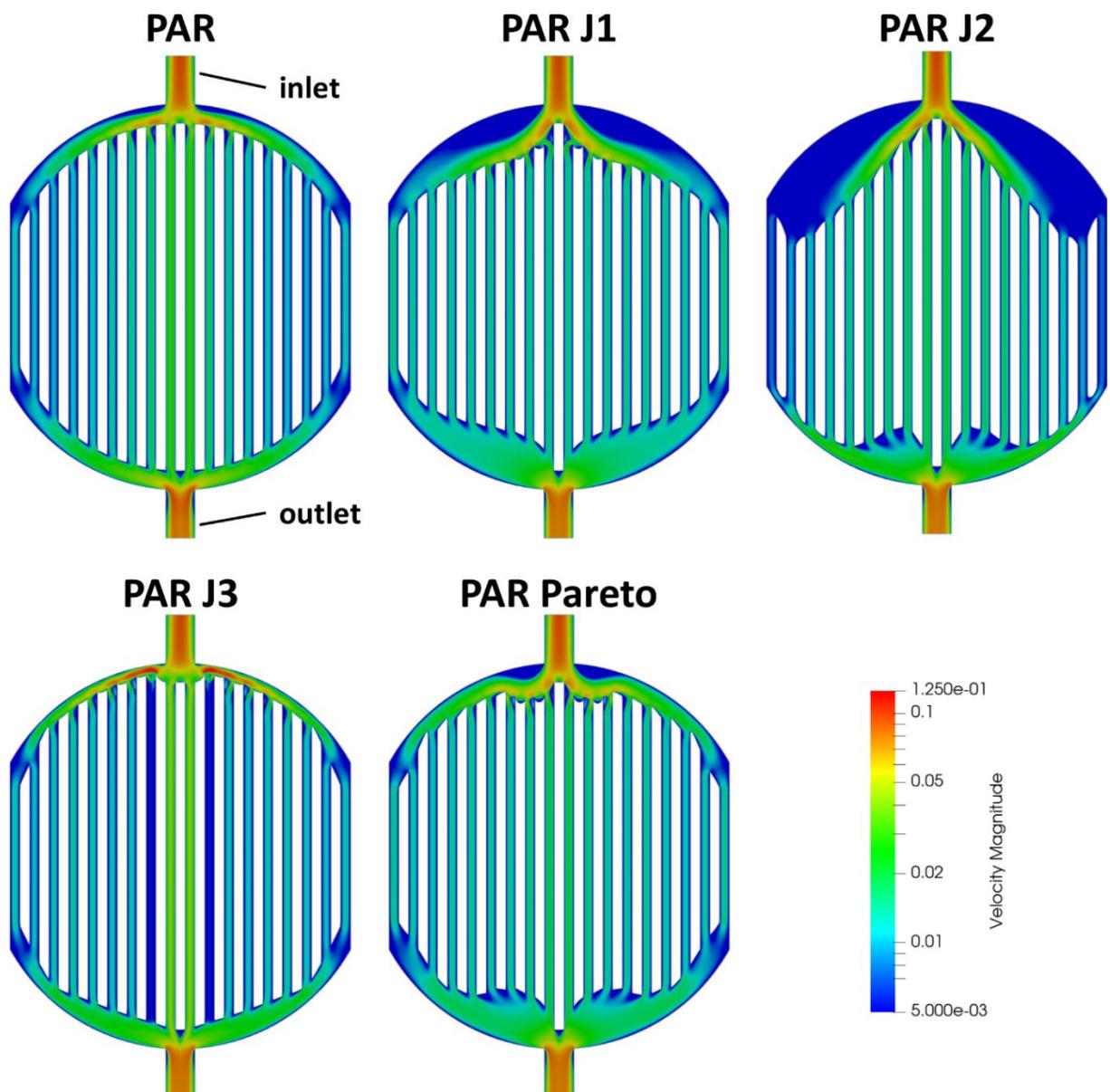

**Figure 4.** Overview of the flow field shapes considered in this paper. The flow velocity is indicated by color using a logarithmic scale. PAR - original shape, PAR J1 - uniform flow, PAR J2 - uniform residence time, PAR J3 - wall shear stress, PAR Pareto - Pareto optimal compromise.

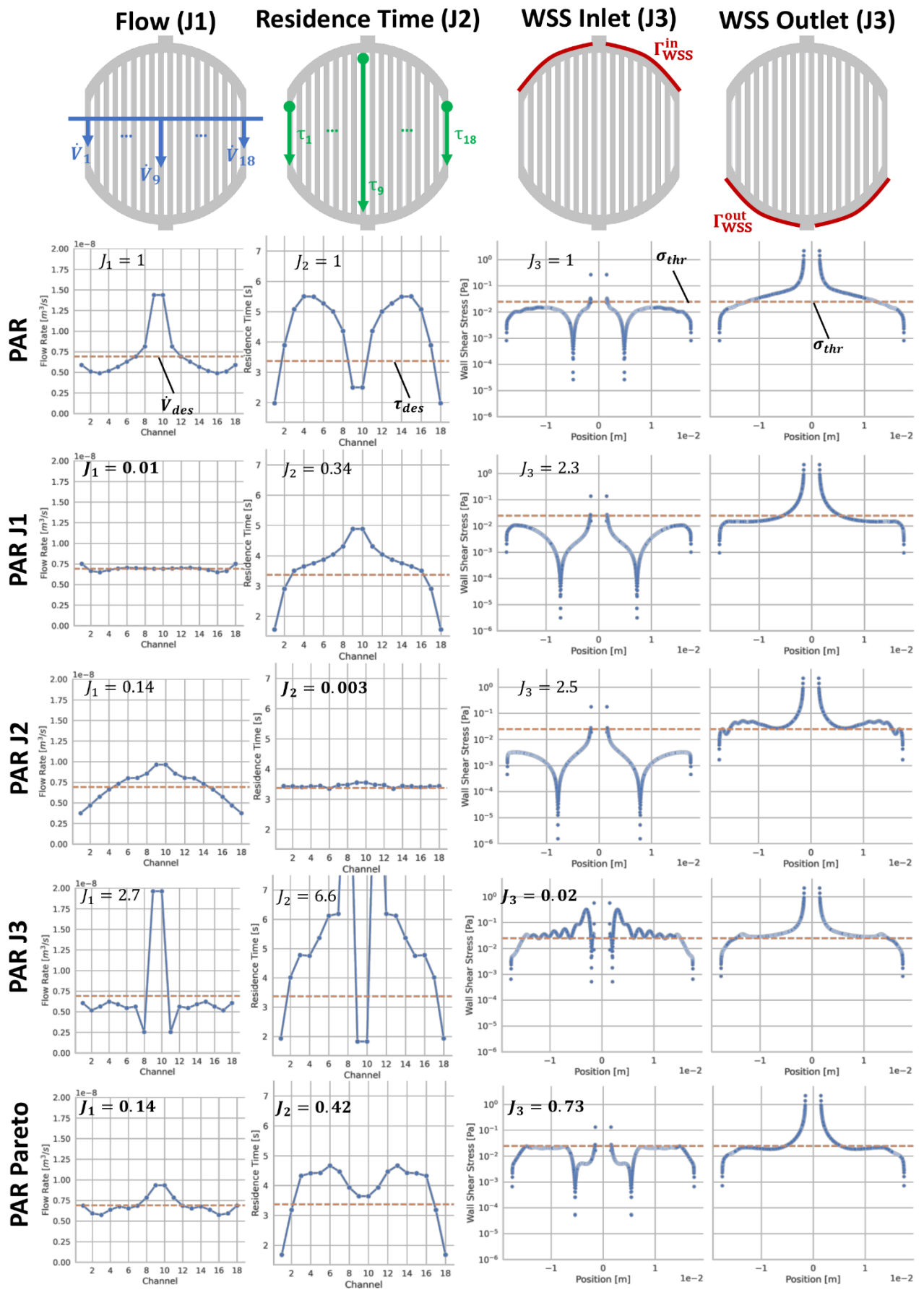

**Figure 5.** Comparison of the flow quantities (velocity, residence time, and wall shear stress) corresponding the optimization criteria $J_1$, $J_2$, and $J_3$ for the five flow field shapes considered.

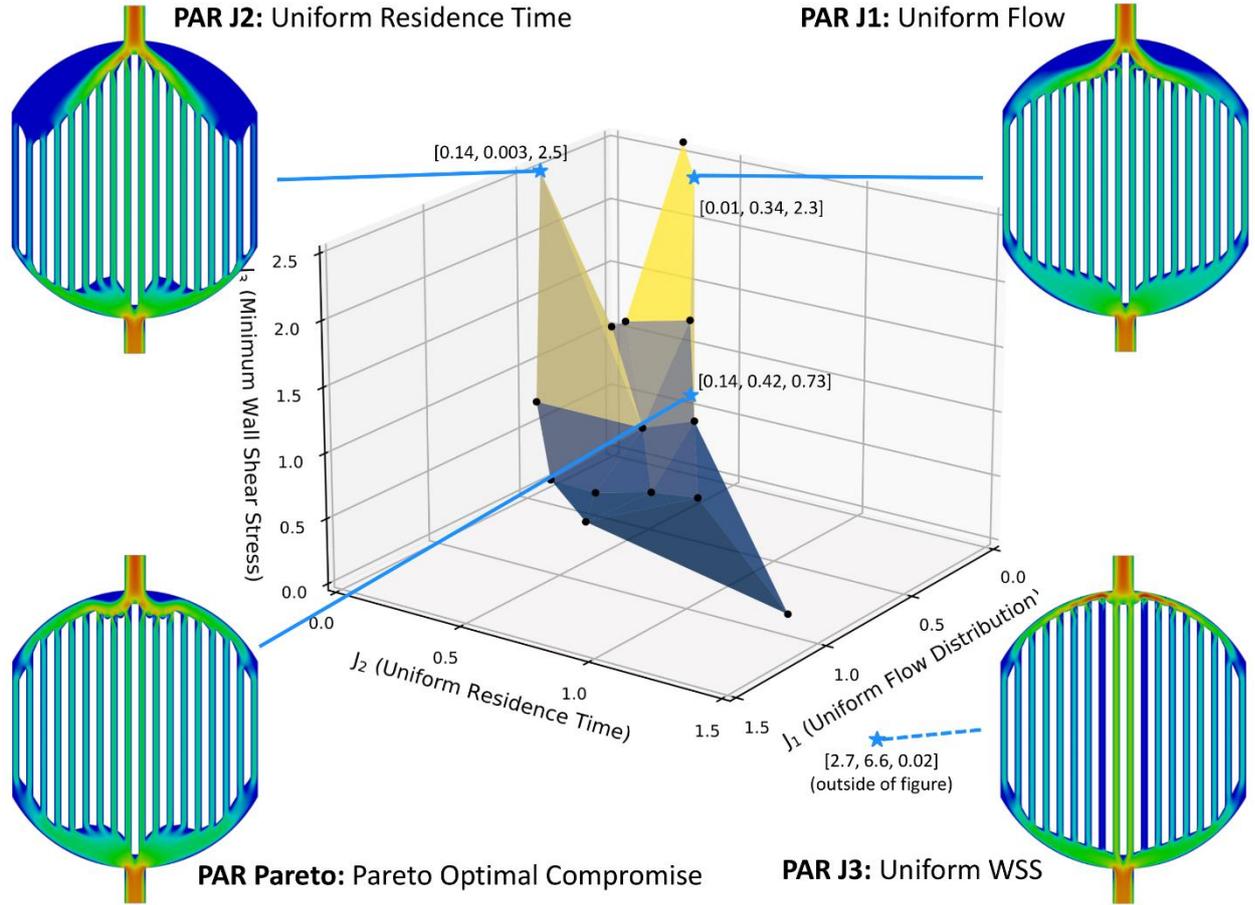

**Figure 6.** Pareto-front and optimal shapes.

### 3.3.1   PAR J1: Flow Field with Uniform Flow
The flow field shape PAR J1 is obtained by solving the single criteria shape optimization problem (11) with cost function $J_1$ (uniform flow). The desired flow rate per channel is chosen to be

$$\dot{V}_{\text{des}} = 6.9 \times 10^{-9} \frac{\text{m}^3}{\text{s}} = \frac{\dot{V}_{in}}{18}, \tag{14}$$

with 18 being the number of channels. The analysis in Figure 5 shows that it is indeed possible to achieve nearly uniform volume flow rates through all channels. This can also be seen from Figure 4: The velocity in the channels is very uniform, which is not the case for the original flow field shape PAR. The inlet and outlet distributor volumes, before and after the channels, respectively, have increased when comparing PAR J1 with the original flow field PAR. The flow pattern in the inlet distributor is different than in the outlet distributor volume due to inertia effects. Especially in the inlet distributor volume of PAR J1 there is a large blue region (see Figure 4) of stagnating flow, which is not desirable in a flow field shape. This stagnation zone can also be seen in the plot of wall shear stress in Figure 5, where the wall shear stress of PAR J1 is significantly lower compared to the original shape PAR. To prevent such a stagnation zone, we later use a combination of cost functions $J_1$ and $J_3$. Before doing so, we first continue to study the influence of each cost function on the geometry, individually.

### 3.3.2   PAR J2: Flow Field with Uniform Residence Time
The shape PAR J2 is obtained by solving the shape optimization problem (11) with cost function $J_2$ (uniform residence time). The desired residence time in the channels is $\tau_{\text{des}} = 3.4$ s. Figure 5 shows a very uniform residence time in all channels. However, the inlet distributor volume is greatly enlarged with a large stagnation area which can also be seen from the low wall shear stress on $\Gamma_{\text{wss}}^{\text{in}}$ (see Figure 5). The flow rate through the outer channels is even lower than for the PAR flow field. The reason is that the outer channels are shorter than the inner channels and a uniform residence time is achieved by

arXiv

reducing the flow rate in the shorter channels. So, if there are channels with different lengths, the criteria $J_1$ and $J_2$ are contradicting.

### 3.3.3 PAR J3: Flow Field with Sufficiently High Wall Shear Stress

The shape PAR J3 is obtained by solving the shape optimization problem (11) with cost function $J_3$. The goal of this cost function is to achieve a sufficiently high wall shear stress on $\Gamma_{wss}$, which should be above the chosen threshold value of $\sigma_{thr} = 0.025$ Pa. The plots in Figure 5 show that the wall shear stress for the PAR J3 shape is indeed mostly above the threshold value. However, the differences in flow rate and residence time have greatly increased. In Figure 4 we observe that particularly the inlet distributor volume decreases slightly to achieve the threshold value for the wall shear stress, confirming our assumptions in Section 2.3.3. Therefore, the wall shear stress criterion can be combined with another criterion to prevent the creation of stagnation zones, which we do in the following.

### 3.3.4 Multi-Criteria Shape Optimization

The Pareto front of (12) is approximated by solving multiple, i.e., 16, scalarized shape optimization problems (13) with different weight vectors to achieve an approximation quality of 0.01 [39] . The approximation of the Pareto front is visualized in Figure 6, where the optimized flow field shapes are shown together with the corresponding point on the Pareto front.

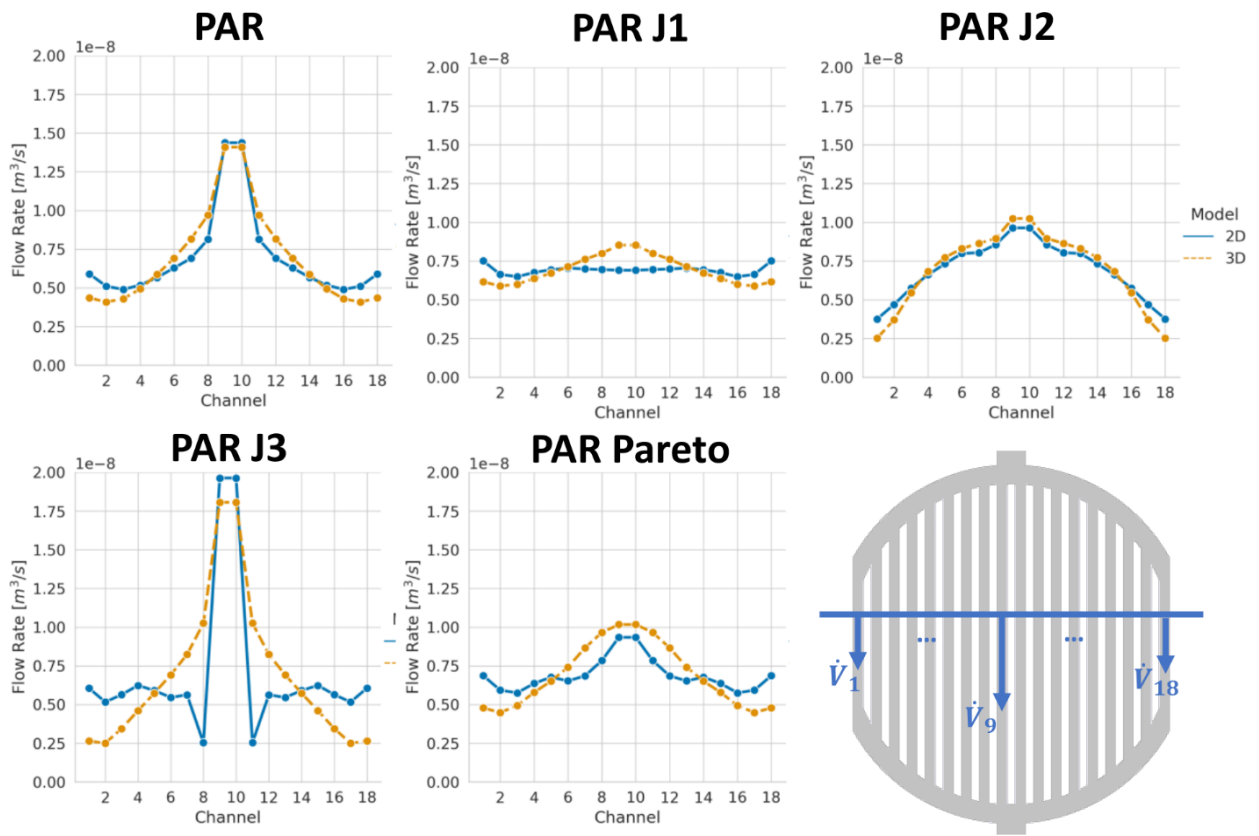

**Figure 8.** Comparison of the flow rate per channel between simplified model (denoted 2D) and detailed model (denoted 3D).

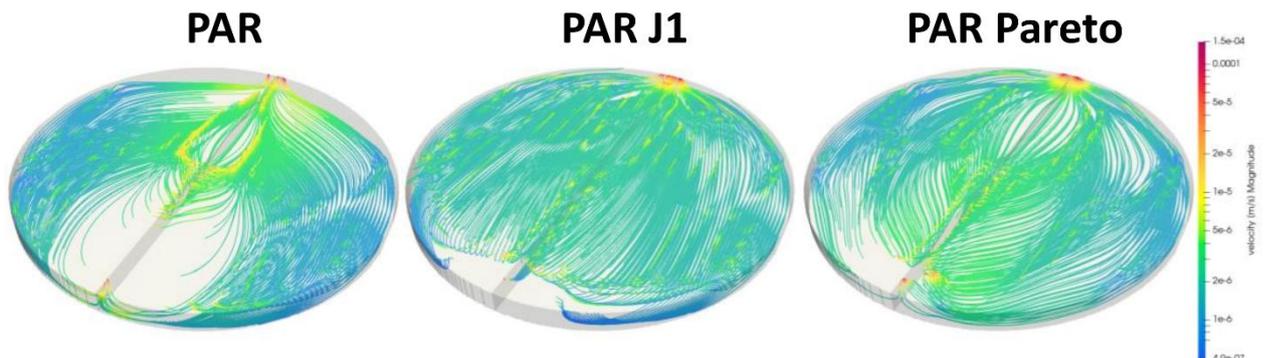

**Figure 7.** Comparison of the flow in the PTL for three different designs.

We have encountered some issues when computing the Pareto front: For some combinations of weights, the shape optimization algorithm did not fully converge. This happened, e.g., when transitioning from cost function $J_1$ to $J_2$, and, hence, the point corresponding to PAR J2 is rather isolated. Usually, the reason for the failed convergence of the algorithm is a decreasing mesh quality. The impact of mesh quality is further discussed, e.g., in [27, 42, 43]. In case of a non-converging scalarized shape-optimization problem (13), the corresponding weight vector provided by the sandwiching algorithm is slightly modified by hand, and then the optimization problem is solved again with the modified weights.

### 3.3.5 PAR Pareto: Selected Pareto Optimal Compromise

To find an optimal compromise between the considered cost functions, we explore the computed Pareto front and investigate the obtained Pareto optimal shapes. We select the point with objective value $J = [0.14, 0.42, 0.73]$ from the Pareto front and re-solve the scalarized shape optimization problem (13) to obtain the corresponding shape which we denote by PAR Pareto. We choose this Pareto solution, i.e., reactor shape, as it seems to show a good uniform flow distribution, due to the cost function $J_1$, but with smaller distributor volumes at the in- and outlet compared to PAR J1, due to the wall shear stress functional $J_3$.

### 3.4 Simulation-Based Benchmarking of Flow Field Shapes

The optimized flow fields investigated in this paper are obtained with the simplified model (3) which does not consider the flow in the PTL. For this reason, we now simulate these flow fields with the detailed simulation model (1) to investigate the influence of the PTL. Results for the comparison of the flow rates per channel are shown in Figure 8. There, some differences in flow rates between the two models are visible. For the flow field PAR J1 the distribution with the more realistic model is not as uniform as for the simplified model, but still good. The reason for the difference in the outer channels is that part of the flow passes through the PTL and can shortcut the way from inlet to outlet. The flow in the middle channels might be higher in some cases since the no-slip condition at the interface is not used here. Better results can be obtained if the model used for shape optimization does also account for the flow in the PTL, which is up for future work. Alternatively, a space mapping approach as demonstrated in [23] can be used.

Finally, Figure 7 shows a comparison of the flow in the PTL between the original flow field PAR and the flow fields PAR J1 and PAR Pareto which have been optimized to make the flow more uniform. The flow velocity is indicated by color on a logarithmic scale. Based on the results shown there, we can conclude that a uniform flow distribution in the flow field also leads to a more uniform flow through the PTL.

## 4 Conclusions

We have shown that our shape optimization approach works well and is a flexible tool for designing novel flow field shapes. As optimality criteria we considered the uniformity of the flow and residence time as well as the wall shear stress. We have seen that methods from multi-criteria optimization are well-suited to select compromises from multiple objective criteria. With the Pareto front at hand, one can easily explore and compare Pareto optimal shapes before making an informed decision. The next step is to experimentally test how the flow fields optimized here affect cell efficiency. For the future, we plan to implement criteria that are more directly related to cell efficiency. For this, the flow and possibly reactions within the PTL must be considered in the optimization problem.

**Acknowledgment**

This work was funded by the Fraunhofer Lighthouse project ShaPID. J. K. and U.-P. A. acknowledge financial support by the cluster of excellence RESOLV, funded by the Deutsche Forschungsgemeinschaft (DFG, German Research Foundation) under Germany's Excellence Strategy – EXC2033-390677874-RESOLV as well as support by the Mercator Research Center Ruhr (MERCUR.Exzellenz, "DIMENSION" Ex-2021-0023 and "KataSign" Ko-2021-0016). D. S. is grateful for financial support by the BMBF (Federal ministry of education and research) within the NanoMatFutur Project "H2Organic" (No. 03XP0421)). K. P. gratefully acknowledges a PhD fellowship granted by the Fonds of the Chemical Industry.

**Symbols used**

| | | |
|---|---|---|
| $Dv$ | [-] | Jacobian of a function $v$ |
| $J$ | [-] | cost functional |
| $K$ | [m²] | permeability |
| $L$ | [m] | length |
| $n$ | [-] | outer unit normal vector on $\Gamma$ |
| $p$ | [Pa] | pressure |
| $u$ | [m/s] | fluid velocity |
| $u_y$ | [m/s] | $y$-component of $u$ |
| $\bar{v}$ | [m/s] | average velocity |
| $\dot{V}$ | [m³/s] | volumetric flow rate |
| $W_\lambda$ | [-] | set of weights for the scalarization of the multi-criteria optimization problem |

**Greek letters**

| | | |
|---|---|---|
| $\Gamma$ | [-] | boundary of the domain $\Omega$ |
| $\lambda$ | [-] | weights for the scalarization of the multi-objective optimization problem |
| $\mu$ | [Pa s] | dynamic viscosity |
| $\rho$ | [kg/m³] | density |
| $\sigma$ | [Pa] | wall shear stress |
| $\tau$ | [s] | residence time |
| $\Omega$ | [-] | computational domain of the system of PDEs |

**Sub- and Superscripts**

| | |
|---|---|
| channels | belonging to the channels of the geometry |
| des | desired |
| eff | effective |
| f | belonging to the fluid part of the domain |
| in | belonging to the inlet |
| int | belonging to the interface between fluid and porous medium |
| out | belonging to the outlet |
| p | belonging to the porous part of the domain |
| sym | belonging to the symmetry boundary |
| thr | threshold |
| wall | belonging to the wall boundary |
| wss | belonging to the part of the boundary where the wall shear stress is considered |

**Abbreviations**

BFGS method – Quasi-Newton method of Broyden, Fletcher, Goldfarb, and Shanno

CFD – computational fluid dynamics

EChH – electrochemical hydrogenation

PAR — name used for the parallel channel geometry considered in this paper

PDE – partial differential equation

PTL – porous transport layer

ShaPID – Name of the Fraunhofer Lighthouse Projekt (Shaping the Future of Green Chemistry by Process Intensification and Digitalization)

ZGE – zero-gap electrolyser

**References**

[1] D. Siegmund, S. Metz, V. Peinecke, T. E. Warner, C. Cremers, A. Grevé, T. Smolinka, D. Segets, U.-P. Apfel, *JACS Au* **2021**, *1 (5)*, 527 – 535. DOI: https://doi.org/10.1021/jacsau.1c00092

[2] S. A. Akhade, N. Singh, O. Y. Gutiérrez, J. Lopez-Ruiz, H. Wang, J. D. Holladay, Y. Liu, A. Karkamkar, R. S. Weber, A. B. Padmaperuma, M.-S. Lee, G. A. Whyatt, M. Elliott, J. E. Holladay, J. L. Male, J. A. Lercher, R. Rousseau, V.-A. Glezakou, *Chemical reviews* **2020**, *120 (20)*, 11370 – 11419. DOI: https://doi.org/10.1021/acs.chemrev.0c00158

[3] Y. Wu, C. Liu, C. Wang, Y. Yu, Y. Shi, B. Zhang, *Nature communications* **2021**, *12 (1)*, 3881. DOI: https://doi.org/10.1038/s41467-021-24059-y

[4] J. Bu, S. Chang, J. Li, S. Yang, W. Ma, Z. Liu, S. An, Y. Wang, Z. Li, J. Zhang, *Nature communications* **2023**, *14 (1)*, 1533. DOI: https://doi.org/10.1038/s41467-023-37251-z

[5] G. Han, G. Li, Y. Sun, *Nat Catal* **2023**, *6 (3)*, 224 – 233. DOI: https://doi.org/10.1038/s41929-023-00923-6

[6] R. S. Sherbo, R. S. Delima, V. A. Chiykowski, B. P. MacLeod, C. P. Berlinguette, *Nat Catal* **2018**, *1 (7)*, 501 – 507. DOI: https://doi.org/10.1038/s41929-018-0083-8

[7] R. S. Sherbo, A. Kurimoto, C. M. Brown, C. P. Berlinguette, *Journal of the American Chemical Society* **2019**, *141 (19)*, 7815 – 7821. DOI: https://doi.org/10.1021/jacs.9b01442

[8] K. Koh, U. Sanyal, M.-S. Lee, G. Cheng, M. Song, V.-A. Glezakou, Y. Liu, D. Li, R. Rousseau, O. Y. Gutiérrez, A. Karkamkar, M. Derewinski, J. A. Lercher, *Angewandte Chemie (International ed. in English)* **2020**, *59 (4)*, 1501 – 1505. DOI: https://doi.org/10.1002/anie.201912241

[9] T. Imada, M. Chiku, E. Higuchi, H. Inoue, *ACS Catal.* **2020**, *10 (22)*, 13718 – 13728. DOI: https://doi.org/10.1021/acscatal.0c04179

[10] K. Pellumbi, L. Wickert, J. T. Kleinhaus, J. Wolf, A. Leonard, D. Tetzlaff, R. Goy, J. A. Medlock, K. Junge Puring, R. Cao, D. Siegmund, U.-P. Apfel, *Chemical science* **2022**, *13 (42)*, 12461 – 12468. DOI: https://doi.org/10.1039/D2SC04647D

arXiv

[11] K. Pellumbi, J. Wolf, S. C. Viswanathan, L. Wickert, M.-A. Kräenbring, J. T. Kleinhaus, K. Junge Puring, F. Özcan, D. Segets, U.-P. Apfel, D. Siegmund, *RSC Sustain.* **2023**, *1 (3)*, 631 – 639. DOI: https://doi.org/10.1039/D3SU00043E

[12] M. Carmo, D. L. Fritz, J. Mergel, D. Stolten, *International Journal of Hydrogen Energy* **2013**, *38 (12)*, 4901 – 4934. DOI: https://doi.org/10.1016/j.ijhydene.2013.01.151

[13] R. Phillips, C. W. Dunnill, *RSC Adv.* **2016**, *6 (102)*, 100643 – 100651. DOI: https://doi.org/10.1039/C6RA22242K

[14] J. D. Egbert, E. C. Thomsen, S. A. O'Neill-Slawecki, D. M. Mans, D. C. Leitch, L. J. Edwards, C. E. Wade, R. S. Weber, *Org. Process Res. Dev.* **2019**, *23 (9)*, 1803 – 1812. DOI: https://doi.org/10.1021/acs.oprd.8b00379

[15] K. Nagasawa, A. Kato, Y. Nishiki, Y. Matsumura, M. Atobe, S. Mitsushima, *Electrochimica Acta* **2017**, *246*, 459 – 465. DOI: https://doi.org/10.1016/j.electacta.2017.06.081

[16] S. Schmidt, C. Ilic, V. Schulz, N. R. Gauger, *AIAA Journal* **2013**, *51 (11)*, 2615 – 2627. DOI: https://doi.org/10.2514/1.J052245

[17] P. Gangl, U. Langer, A. Laurain, H. Meftahi, K. Sturm, *SIAM J. Sci. Comput.* **2015**, *37 (6)*, B1002-B1025. DOI: https://doi.org/10.1137/15100477X

[18] G. Leugering, A. A. Novotny, G. P. Menzala, J. Sokołowski, *Math. Meth. Appl. Sci.* **2010**, *33 (17)*, 2118 – 2131. DOI: https://doi.org/10.1002/mma.1324

[19] R. Hohmann, C. Leithäuser, *SIAM J. Sci. Comput.* **2019**, *41 (4)*, B625-B648. DOI: https://doi.org/10.1137/18M1225847

[20] C. Leithäuser, R. Pinnau, R. Feßler, *Optim Eng* **2018**, *19 (3)*, 733 – 764. DOI: https://doi.org/10.1007/s11081-018-9396-3

[21] S. Blauth, C. Leithäuser, R. Pinnau, *Z Angew Math Mech* **2021**, *101 (4)*. DOI: https://doi.org/10.1002/zamm.202000166

[22] S. Blauth, C. Leithäuser, R. Pinnau, *Journal of Mathematical Analysis and Applications* **2020**, *492 (2)*, 124476. DOI: https://doi.org/10.1016/j.jmaa.2020.124476

[23] S. Blauth, *SIAM J. Optim.* **2023**, *33 (3)*, 1707 – 1733. DOI: https://doi.org/10.1137/22M1515665

[24] S. Blauth, *SoftwareX* **2021**, *13*, 100646. DOI: https://doi.org/10.1016/j.softx.2020.100646

[25] V. H. Schulz, M. Siebenborn, K. Welker, *SIAM J. Optim.* **2016**, *26 (4)*, 2800 – 2819. DOI: https://doi.org/10.1137/15M1029369

[26] S. Blauth, *SIAM J. Optim.* **2021**, *31 (3)*, 1658 – 1689. DOI: https://doi.org/10.1137/20M1367738

[27] P. M. Müller, N. Kühl, M. Siebenborn, K. Deckelnick, M. Hinze, T. Rung, *Struct Multidisc Optim* **2021**, *64 (6)*, 3489 – 3503. DOI: https://doi.org/10.1007/s00158-021-03030-x

[28] M. Nemec, D. W. Zingg, T. H. Pulliam, *AIAA Journal* **2004**, *42 (6)*, 1057 – 1065. DOI: https://doi.org/10.2514/1.10415

[29] R. Hilbert, G. Janiga, R. Baron, D. Thévenin, *International Journal of Heat and Mass Transfer* **2006**, *49 (15-16)*, 2567 – 2577. DOI: https://doi.org/10.1016/j.ijheatmasstransfer.2005.12.015

[30] P. Di Barba, M. E. Mognaschi, *IEEE Trans. Magn.* **2009**, *45 (3)*, 1482 – 1485. DOI: https://doi.org/10.1109/TMAG.2009.2012685

[31] I. Vasilopoulos, V. G. Asouti, K. C. Giannakoglou, M. Meyer, *Engineering with Computers* **2021**, *37 (1)*, 449 – 459. DOI: https://doi.org/10.1007/s00366-019-00832-y

[32] O. Iliev, R. Kirsch, Z. Lakdawala, S. Rief, K. Steiner, in.

[33] A. Parasyris, C. Brady, D. B. Das, M. Discacciati, *AMST* **2019**, *23 (3)*. DOI: https://doi.org/10.11113/amst.v23n3.158

[34] https://www.itwm.fraunhofer.de/en/departments/sms/products-services/filtest-filter-element-simulation-toolbox.html (Accessed on September 12, 2023).

[35] M. Hinze, R. Pinnau, M. Ulbrich, S. Ulbrich, *Optimization with PDE Constraints*, Vol. 23, Springer Netherlands, Dordrecht **2009**.

[36] M. Alnaes, J. Blechta, J. Hake, A. Johansson, B. Kehlet, A. Logg, C. Richardson, J. Ring, M. E. Rognes, G. N. Wells, *Archive of numerical software* **2015**, *3 (100)*, 9 – 23.

[37] Matthias Ehrgott, *Multicriteria Optimization*, Springer-Verlag, Berlin/Heidelberg **2005**.

[38] *Multicriteria optimization in intensity modulated radiotherapy planning*, Berichte des Fraunhofer-Instituts für Techno- und Wirtschaftsmathematik (ITWM Report), Vol. 77, Fraunhofer-Institut für Techno- und Wirtschaftsmathematik, Fraunhofer (ITWM), Kaiserslautern **2005**.

[39] M. Bortz, J. Burger, N. Asprion, S. Blagov, R. Böttcher, U. Nowak, A. Scheithauer, R. Welke, K.-H. Küfer, H. Hasse, *Computers & Chemical Engineering* **2014**, *60*, 354 – 363. DOI: https://doi.org/10.1016/j.compchemeng.2013.09.015

[40] D. Nowak, K.-H. Küfer, in press.

[41] R. Banerjee, N. Bevilacqua, L. Eifert, R. Zeis, *Journal of Energy Storage* **2019**, *21*, 163 – 171. DOI: https://doi.org/10.1016/j.est.2018.11.014

[42] R. P. Dwight, in *Computational Fluid Dynamics 2006* (Eds: H. Deconinck, E. Dick), Springer Berlin Heidelberg. Berlin, Heidelberg **2009**.

[43] T. Etling, R. Herzog, E. Loayza, G. Wachsmuth, *SIAM J. Sci. Comput.* **2020**, *42 (2)*, A1200-A1225. DOI: https://doi.org/10.1137/19M1241465